\documentclass[12pt]{amsart}

\newtheorem{theorem}{Theorem}[section]
\newtheorem{lemma}[theorem]{Lemma}
\newtheorem{claim}[theorem]{Claim}

\newtheorem{proposition}[theorem]{Proposition}

\theoremstyle{definition}

\newtheorem{openpb}[theorem]{Open Problem}

\theoremstyle{remark}
\newtheorem{remark}[theorem]{Remark}

\numberwithin{equation}{section}

\newcommand{\Hh}{\mathbb{H}}
\newcommand{\D}{\mathbb{D}}
\newcommand{\C}{\mathbb{C}}
\newcommand{\N}{\mathbb{N}}
\newcommand{\B}{\mathbb{B}}

\newcommand{\R}{\mathbb{R}}

\newcommand{\De}{\Delta }

\newcommand{\si}{\sigma }

\newcommand{\ga}{\gamma }
\newcommand{\Ga}{\Gamma }

\newcommand{\Om}{\Omega }

\newcommand{\rea}{\operatorname{Re}}
\newcommand{\ima}{\operatorname{Im}}
\newcommand{\Arg}{\operatorname{Arg}}

\newcommand{\ntlim}{\operatorname{\hbox{K-}lim}}

\newcommand{\bd}[1]{\partial #1}

\begin{document}
\baselineskip=18pt
%\large

\title[Valiron's Theorem]{On Valiron's Theorem}

\author{Filippo Bracci}
\address{Dipartimento di Matematica, Universit\`a di Roma ``Tor
Vergata'', Via della Ricerca Scientifica 1, 00133 Roma, Italy.}
\email{fbracci@mat.uniroma2.it}

\author{Pietro Poggi-Corradini}
\address{Department of Mathematics, Cardwell Hall, Kansas State University,
Manhattan, KS 66506, USA.} \email{pietro@math.ksu.edu}

%\dedicatory{PRELIMINARY VERSION}
\date{\today}

\begin{abstract}
This is a survey on Valiron's Theorem about the convergence
properties of orbits of analytic self-maps of the disk of
hyperbolic type and related questions in one and several
variables.
\end{abstract}

\maketitle

\section{Introduction}\label{sec:intro}

Let $\D=\{z\in \C: |z|<1 \}$ and let $\phi$ be  an analytic
function defined on $\D$. If $|\phi (z)|<1$ for $|z|<1$, then
$\phi$ is a self-map of the disk $\D$ and one can iterate by
letting $\phi_{n}=\phi \circ \cdots\phi$, $n$ times. The natural
question that arises is {\em given a point $z_{0}\in \D$, what can
be said about its orbit $z_{n}=\phi_n(z_0)$, as $n=1,2,3,\dots $?}
In this survey we will describe a theorem of Valiron which relates
to this question and describe the multidimensional setting.

\subsection{Schwarz's Lemma}\label{ssec:sl}
One of the very first results one encounters in function theory is
Schwarz's Lemma, which can be proved using the maximum principle.
\begin{lemma}[Schwarz's Lemma]\label{lem:sl}
Suppose $\phi$ is an analytic self-map of $\D$. If, moreover, $\phi
(0)=0$ then
\begin{enumerate}
\item $|\phi (z)|\leq |z| $ for all $z\in \D$.
\item $|\phi (z_{0})|=|z_{0}|$ for some $z_{0}\neq 0$ if and only if
$\phi$ is a rotation.
\item $|\phi^{\prime} (0)|\leq 1$ and $|\phi^{\prime} (0)|= 1$ if and
only if $\phi$ is a rotation.
\end{enumerate}
\end{lemma}
The proof is based on the fact that the function $\phi (z)/z$ is
analytic and bounded by $1$.

Geometrically, Schwarz's Lemma says that for every $0<r<1$:
\[
\phi (r\D)\subset r\D,
\]
and from the proof one deduces more precisely that,
except for rotations, for every $0<r_{0}<1$ there
exists $0<s_{0}<1$ such that for $0<r<r_{0}$,
\begin{equation}\label{eq:contraction}
\phi (r\D)\subset s_{0}r\D.
\end{equation}

The maximum principle and Schwarz's Lemma can be used to show that the
automorphisms of $\D$ are of the form
\[
\gamma (z)=c\frac{z-a}{1-\bar{a}z}
\]
for some constants $|c|=1$ and $|a|<1$. All the automorphisms $\ga $
send the point $a\in\D$ to $0$, and we write $\ga_{a}$ when the
constant $c$ equals $1$.
Using these automorphisms one can transfer the system of disks $r\D$
($0<r<1$) around any given point  $a\in\D$:
\[
\Delta(a,r)=\ga_{a}^{-1} (r\D).
\]
These are called pseudo-hyperbolic disks of radius $r$ at $a$.
Because linear fractional transformation map circles to circles,
$\Delta (a,r)$ is an Euclidean disk, but $a$ is not the Euclidean
center (actually $a$ is further from the origin). With these
notations, a simple use of the $\ga_{z}$ and $\ga_{\phi (z)}$ shows
that given any analytic self-map $\phi$ of $\D$, and for any $z\in\D$,
we always have,
\begin{equation}\label{eq:ifsl}
\phi (\Delta (z,r))\subset \Delta (\phi(z),r)
\end{equation}
for all $0<r<1$. Moreover, by continuity and compactness,
given a compact set $K=\{|z|\leq t \}$ for some $0<t<1$, and given a
radius $0<r_{0}<1$  there exists a constant $s_{0}<1$ so that uniformly
for $z\in K$ and for $0<r<r_{0}$:
\begin{equation}\label{eq:contraction2}
\phi (\Delta (z,r))\subset \Delta (\phi(z),s_{0}r).
\end{equation}

This can also be worded in terms of the pseudo-hyperbolic distance
\[
d (z,w)=|\ga_{z} (w)|=\left|\frac{z-w}{1-\bar{w}z}\right|\qquad
\mbox{ for } z,w\in \D.
\]
Although we call it distance, $d (z,w)$ does not satisfy the triangle
inequality, yet it almost does for small distances because of the formula:
\[
d (z,w)\leq \frac{d (z,\zeta)+d (\zeta ,w)}{1+d (z,\zeta)d (\zeta
,w)}.
\]
An actual distance is obtained by letting
\[
\rho (z,w)=\log\frac{1+d (z,w)}{1-d (z,w)}.
\]
This is the hyperbolic distance of $\D$.

\subsection{One fixed point in $\D$}\label{ssec:onefp}
If a self-map of the disk fixes two points, conjugating it  using an
automorphism and using part 2. of Schwarz's Lemma \ref{lem:sl}, one
proves that it is actually the identity map. Thus every other self-map
of the disk can fix at most one point in $\D$.

If a self-map fixes exactly one point in $\D$ it is called of {\em
elliptic type}. In this case, the map can be conjugated by
automorphisms so that the fixed point is the origin, hence the power
series expansion there is:
\[
\phi (z)=\lambda z + O (z^{2})
\]
where $\lambda =\phi^{\prime} (0)\in \overline{\D}$. Three subcases
arise: when $|\lambda |=1$, the map $\phi$ is a rotation; if
$0<|\lambda |<1$ the fixed point is called {\em attractive}, if
$\lambda =0$ it is {\em superattractive}. The behavior of single
orbits is well understood in all these cases.

For instance, in the attractive case, it is clear from
(\ref{eq:contraction}) that every orbit tends to zero. Moreover,
by K\oe nigs Theorem (which is proved using property
(\ref{eq:contraction2})), there is a one-to-one analytic map
$\sigma$ defined near $0$ with $\si (0)=0$ and   $\si^{\prime}
(0)=1$ such that
\begin{equation}\label{eq:koenigs}
\si \circ \phi  (z)=\lambda \si(z)
\end{equation}
near $0$, i.e one can change coordinates holomorphically so that
$\phi$ becomes linear, and because of the condition on
$\si^{\prime} (0)$ being equal to $1$, the orbits $z_{n}=\phi_{n}
(z_{0})$ asymptotically approach the corresponding orbit
$\lambda^{n}\si (z_{0})$. It is worth notice that $\sigma$ can be
extended (not univalently in general) to all $\D$ in such a way
that \eqref{eq:koenigs} is still holding.

In the superattractive case,  the orbits tend fast to the origin.
Even in this case it is possible to perform an holomorphic  change
of variables  near the origin in such a way that $\phi$ assumes a
simpler form. Namely, if $\phi(z)=O(z^k)$ then by B\"ottcher's
Theorem there exists a one-to-one analytic map $\sigma$ defined
near $0$ so that $\sigma(0)=0$, $\si'(0)=1$ and
\begin{equation}\label{eq:bottcher}
\si \circ \phi (z)= \si(z)^k
\end{equation}
near $0$. In this case however the map $\sigma$ cannot be in
general well defined on all $\D$. For all these matters we refer
the interested reader to~\cite{cg}.

\subsection{No fixed points in $\D$}\label{ssec:nofp}
Assume that the self-map $\phi$ fixes no point in $\D$. Then either $\phi$
is an automorphism of $\D$ in which case it is an isometry for the
hyperbolic distance, or, by Schwarz's Lemma,
$\phi$ is a strict contraction, i.e.
\[
d (\phi(z),\phi(w))<d (z,w).
\]
for all $z,w \in\D$. If $\phi$ is an automorphism then it can be
conjugated to one of two maps: either multiplication by $T>1$ on the
upper half-plane $\Hh =\{\ima z>0 \}$ ({\em hyperbolic automorphism})
or translation by $b>0$ on $\Hh$ ({\em parabolic automorphism}).

We will see that $\phi$
can also be classified as
hyperbolic or parabolic when it is not an automorphism.
However, even though self-maps of the disk with no fixed points
do try to imitate
the behavior of the automorphisms in the long run, this is only true to
varying degrees and the situation is much more complicated, especially
in the parabolic case.
The main topic of this survey is to describe self-maps $\phi$ of
hyperbolic type.

We will proceed in stages. The first claim is that given a
self-map of the disk there exists a point $\zeta \in\bd{\D}$ such
that every orbit of $\phi$ converges to $\zeta$. This allows one
to change variables to the upper half-plane and send $\zeta$ to
infinity. Computation usually become easier in this formulation,
although it might still be useful to work in both models in view
of the possible extentions to several complex variables. The point
$\zeta$ is the famous {\em Denjoy-Wolff point} of the map $\phi$.
The second claim is that, like in the case of automorphisms, if
$\phi$ (not an automorphism) is in the upper half-plane model with
Denjoy-Wolff point at infinity, then either $\ima \phi (z)> \ima
z$ for all $z\in \Hh$ (parabolic case), or there is $T>1$ such
that $\ima \phi (z)> T \ima z$ for all $z\in \Hh$ (hyperbolic
case). Notice that if $\phi$ is not an automorphism then the
previous inequalities are strict at {\em every} point. Indeed, by
a generalization of Schwarz's Lemma known as Julia's Lemma (see
below), if there is equality at some point then there is equality
everywhere and $\phi$ is an automorphism of either parabolic type
(in the first) or hyperbolic type (in the second).

 In more geometric terms,
letting $H (s)=\{\ima z>s \}$, $\phi (H (s))\subset H (Ts)$ for
some $T\geq 1$. The half-planes $H (s)$ are called {\em horodisks}
because in the disk model they correspond to Euclidean disks
tangent to $\bd{\D}$ at the Denjoy-Wolff point.

We first observe that if $\phi$
is a self-map of $\D$ with no fixed points and $\phi$ is not an
automorphism, then no iterate of $\phi$ can have fixed points in $\D$ either.
In fact, suppose that $\phi_N (z_{0})=z_{0}$ for some $N\geq 2$ and
some $z_{0}\in \D$. The the orbit of $z_{0}$ is periodic of period $N$
and so are the steps $d_{n}=d (z_{n},z_{n+1})$, which contradicts the
fact that $d_{n}$ is a strictly decreasing sequence by Schwarz's Lemma.

This can be used to show that any orbit $z_{n}$ cannot accumulate
anywhere in $\D$, i.e. must eventually escape any given compact set.
In fact, suppose that a subsequence $z_{n_{k}}$ tends to $p\in \D$.
Find $0<t<1$ so that $|p|<t$ and let $K=\{|z|\leq t \}$, also let
$s<1$ and $0<r_{0}<1$ be given as in (\ref{eq:contraction2}).
Eliminating finitely many terms, we can assume that $z_{n_{k}}\in K$
for all $k$.
Choose a radius $0<r_{0}<1$ close enough to $1$ so that the
pseudo-hyperbolic disk $D=\Delta (z_{n_{1}},r_{0})$ contains $K$.
By (\ref{eq:contraction2}) we have
\[
\phi_{n_{k}} (D)\subset \Delta (z_{n_{k}},s_{0}^{k}r_{0}),
\]
and since $s_{0}^{k}$ tends to zero and $z_{n_{k}}$ tends to $p$, we
must have for large enough $k_{0}$ that $\phi_{n_{k_{0}}}
(\overline{D})\subset D$. This implies that $\phi_{n_{k_{0}}}$ has a
fixed point in $D$ but we have ruled out fixed points for the
iterates of $\phi$.

The next step is to show that given an orbit $z_{n}$, not only
$|z_{n}|$ tends to one but actually $z_{n}$ tends to some $\zeta \in
\bd{\D}$. For this we need a boundary consequence of Schwarz's Lemma
known as Julia's Lemma.

\subsection{Julia's Lemma}\label{ssec:jl}
We will present a simplified version of Julia's Lemma which is more
suitable to our purposes. First we use the Poisson kernel at $\zeta\in\bd\D$
to define the horodisks at $\zeta$:
\begin{equation}\label{eq:horodisks}
H (t)=\{z\in\D: \frac{1-|z|^{2}}{|\zeta -z|^{2}}>\frac{1}{t}\}.
\end{equation}
Note that $H (t)$ is decreasing as $t\downarrow 0$ and
$\cap_{t\downarrow 0}H (t)=\emptyset$ while $\cup_{t\uparrow \infty}H
(t)=\D$.
\begin{lemma}\label{lem:jl}
Let $\phi$ be an analytic  self-map of $\D$ and let $p_{k}\in \D$ be a
sequence of points tending to $\zeta \in \bd{\D}$. If $\phi (p_{k})$
also tends to $\zeta$ and the ratio
\[
\frac{1-|\phi (p_{k})|}{1-|p_{k}|}\longrightarrow d>0,
\]
as $k\rightarrow \infty$, then for all $t>0$
\begin{equation}\label{eq:jl}
\phi (H (t))\subset H (dt).
\end{equation}
\end{lemma}
The general version of Julia's Lemma allows for $\phi(p_n)$ to be
tending to some other boundary point $\eta\in \bd\D$ and also does not
assume $d>0$, but deduces it.
The proof of this lemma is obtained by applying
Schwarz's Lemma in the form (\ref{eq:ifsl}) to hyperbolic disks
centered at $p_{n}$ of larger and larger radius so that these disks
tend to the horodisk $H (t)$. Note also that while (\ref{eq:ifsl})
contracts the hyperbolic radius,
when $d>1$ equation (\ref{eq:jl})  only requires
for a smaller horodisk to be mapped into a larger one.

Now consider an orbit $z_{n}$. We have seen above that $|z_{n}|$ tends
to $1$. Choose a subsequence $z_{n_{k}}$ such that
\[
|z_{n_{k}+1}|=|\phi (z_{n_{k}})|\geq |z_{n_{k}|}
\]
and further assume that $z_{n_{k}}$ tends to some point $\zeta \in
\bd\D$. Since $d (z_{n_{k}},\phi (z_{n_{k}}))\leq d
(z_{0},z_{1})$, $\phi (z_{n_{k}})$ also tends to $\zeta$. Hence,
we can apply Julia's Lemma, with $p_{k}=z_{n_{k}}$ and with $d\leq
1$, to find that $\phi (H (t))\subset H (t)$. This immediately
implies that the whole orbit $z_{n}$ must tend to $\zeta$.
Moreover if we let
\[
\alpha =\liminf_{z\rightarrow \zeta}\frac{1-|\phi (z)|}{1-|z|}
\]
then by Julia's Lemma $\phi (H (t))\subset H (\alpha t)$. We call $\alpha$
 the {\em  coefficient of dilatation} of $\phi$ at its Denjoy-Wolff
point. It follows from what we said so far that $\alpha \leq 1$,
and the map $\phi$ is said to be of {\em hyperbolic type} if
$\alpha <1$, while it is of {\em parabolic type} if $\alpha =1$.
It can also be shown that $\alpha >0$ always.

As we mentioned above the terminology parabolic vs. hyperbolic is used
because one wishes to show that these maps tend to imitate the
corresponding parabolic vs. hyperbolic automorphisms. However, this is
not always the case, especially in the parabolic case. What happens in
the hyperbolic case is the content of this survey.

\section{Self-maps of the disk of hyperbolic type}\label{sec:hyp}

The hyperbolic automorphisms in the upper half-plane
model are easy to describe. They are of the form
\[
\tau (z)= Az +b
\]
with $A>1$ and $b\in \R$. The only two fixed points for $\tau$ are
infinity and $-b/ (A-1)$. The hyperbolic geodesic $L=\{\rea z=-b/
(A-1);\ \ima z>0 \}$ is invariant ($L$ is also known as the axis of
the automorphism) and invariant is also every half line
originating from the fixed point $-b/ (A-1)$ which lays in $\Hh$.
It is clear then that for every orbit $z_{n}$ of $\tau$, the following
three properties hold:
(1) the ratios $z_{n+1}/z_{n}$ tend to $A$;
(2) the sequence $\Arg
z_{n}$ has (a) a limit in $(0,\pi)$ and a computation shows that it equals
$\Arg (z_{0}-b/ (A-1))$, hence it is a harmonic function of $z_{0}$ and
(b) by varying $z_{0}$ this limit takes every value in $(0,\pi)$;
(3) the sequence $z_{n}/A^{n}$ tends $z_{0}-b/ (A-1)$.

A quicker way to describe this dynamic would have been to notice that $\tau$
can be conjugated via a translation to the map
$z\mapsto Az$.
Back in the disk model the axis $L$ is an arc of circle orthogonal to
$\bd{\D}$, intersecting $\bd{\D}$ at $1$ and at some other fixed point
$p\in \bd{\D}\setminus \{1 \}$. All arcs of circle interesecting
$\overline{\D}$ in $1$ and $p$ are invariant for the automorphism
and the three properties above become:
(1) the ratios $(1-z_{n+1})/ (1-z_{n})$ tend to $\alpha $;
(2) the sequence $\Arg (1-z_{n})$ has (a) a limit in $(-\pi/2,\pi/2)$
which is a harmonic function of $z_{0}$ and (b) this limit takes every
value in $(-\pi/2,\pi/2)$;
(3) the sequence $(1-z_{n})/\alpha ^{n}$ tends to a limit.

Assume now that $\phi$ is a self-map of the disk with Denjoy-Wolff point
at $1$ (without loss of generality) and coefficient of dilatation
$\alpha <1$. Or, equivalently, assume that $\Phi$ is a self-map of the upper
half-plane $\Hh$ and $\Phi(z)=Az+p (z)$ with $\ima p (z)\geq 0$ and
$A=1/\alpha= \inf_{z\in\Hh}\frac{\ima \Phi (z)}{\ima z}>1$.
It is natural to ask if the three properties of
hyperbolic automorphisms above are also shared by the orbits of
$\Phi$.
Valiron shows that properties (1) and (2) (a) still hold, see \cite{va1}
or Chapter VI of \cite{va2} (he doesn't seem to have considered
property 2 (b)). Next we
present a slighty different proof of his result.

\subsection{
%\noindent
{\bf Property (1):} {\em Given an orbit $z_{n}$ of $\Phi$, the
ratio $\frac{z_{n+1}}{z_{n}}$ tends to $A$}.
}\label{ssec:jc}

\noindent This property
is intimately connected with the Julia-Carath\'eodory
Theorem. We state it somewhat reworded, in the upper half-plane model.
\begin{theorem}[Julia-Carath\'eodory]\label{thm:jc}
Let $\Phi$ be an analytic self map of $\Hh$. Let
\[
A=\inf_{z\in\Hh}\frac{\ima \Phi (z)}{\ima z}.
\]
Then,
\begin{equation}\label{eq:jc}
\ntlim_{z\rightarrow \infty}\frac{\Phi(z)}{z}=A.
\end{equation}
\end{theorem}
For a proof see \cite{sh} p.~66-69, which, as one might guess, is based
on Schwarz's Lemma.
By $\ntlim_{z\rightarrow \infty}$, ``non-tangential limit'', we mean
that $z$ tends to infinity
in such a way that $|\Arg z-\pi /2|<\pi /2-\delta$ for any given $\delta >0$

In particular, when $\Phi$ is of hyperbolic type then
(\ref{eq:jc}) holds. Yet one cannot immediately deduce from it
property (1) for the orbits of $\Phi$ since, in principle, $z_n$
might tend to $1$ tangentially.

\begin{lemma}\label{lem:ntapproach}
Let $\Phi$ be a hyperbolic holomorphic self map of $\Hh$. Then any
given orbit $z_{n}$ satisfies
\[
|\Arg z_{n}-\pi /2|<\pi /2-\delta
\]
for some fixed $\delta>0$ depending only on $z_{0}$.
\end{lemma}
\begin{proof}
Schwarz's Lemma imposes that
$z_{n+1}$ belongs to the pseudo-hyperbolic disk $\Delta$ centered at
$z_{n}$ of radius $d_{0}=d (z_{0},z_{1})$, and the hyperbolic type imposes
that $\ima z_{n+1}\geq A\ima z_{n}$ for some $A>1$. So $z_{n+1}$ is
forced to land in the intersection (never empty!) of $\Delta$ with the
half-plane $\{\ima z\geq A\ima z_{n} \}$. Applying a dilation $1/\ima z_{n}$ to
this picture we see that $z_{n+1}/\ima z_{n}$ belongs to the
intersection of a pseudo-hyperbolic disk of radius $d_{0}$, centered
at some point with imaginary part equal to $1$, and the half-plane
$\{\ima z\geq A \}$. From this we deduce that
\[
|\Arg (z_{n+1}-z_{n}) -\pi /2|\leq \pi/2-\delta_{0}
\]
for some $\delta_{0} >0$ which  depends only on $z_{0}$. Now
consider a sector $S (\delta ) =\{|\Arg z-\pi /2|\leq \pi
/2-\delta \}$ and let $R$ be the union of all the sectors $z+S
(\delta_{0} )$ as $z$ describes the segment
$[-\overline{z_{0}},z_{0}]$. It is clear that the orbit $z_{n}$
never leaves the region $R$, and that $R$ is contained in a larger
sector $S (\delta_{1})$ with $0<\delta_{1}<\delta_{0}$.
\end{proof}
Now that we know that every orbit stays confined in a non-tangential
approach region, we can apply Julia-Carath\'eodory's
theorem and obtain property (1) that $z_{n+1}/z_{n}$ always tends to $A$.

\subsection{
%\noindent
{\bf Property (2) (a):} {\em Given an orbit $z_{n}$ of $\Phi$, the
limit $\Arg z_{n}$ exists and is a harmonic function of $z_{0}$}.
}\label{ssec:arg}

\noindent Observe first that $\Arg z_{n}=\Arg \Phi_{n}
(z_{0})$ is a bounded harmonic function in $z_{0}$, so once the
existence of the limit is established, harmonicity follows by
Harnack's principle. We write $z_{n}=x_{n}+iy_{n}$.
Property (1) can be written as
$z_{n+1}=Az_{n}+o(1)z_{n}$, thus dividing by $y_{n}$ we get
\[
\frac{z_{n+1}}{y_{n}}=A\frac{z_{n}}{y_{n}}+o
(1)\frac{z_{n}}{y_{n}}.
\]
However, Lemma \ref{lem:ntapproach} implies that
$z_{n}/y_{n}=x_{n}/y_{n}+i$ is bounded away from $0$ and
$\infty$.
So, taking the imaginary part of both sides, we obtain
\begin{equation}\label{eq:yn}
\frac{y_{n+1}}{y_{n}}=A+o (1).
\end{equation}
Consider the automorphism of $\Hh$ that sends $z_{n}$ back to $i$, i.e.
\begin{equation}\label{eq:taun}
\tau_{n} (z)=\frac{z-x_{n}}{y_{n}}.
\end{equation}
Then set
\begin{equation}\label{eq:qn}
q_{n}=\tau_{n}
(z_{n+1})=\frac{x_{n+1}-x_{n}}{y_{n}}+i\frac{y_{n+1}}{y_{n}}.
\end{equation}
It follows immediately from (\ref{eq:yn}) that $\ima q_{n}\rightarrow
A$ as $n$ tends to infinity. Also,
by conformal invariance, the sequence
\[
d (i,q_{n})=d (z_{n},z_{n+1})=d_{n}\geq \frac{1-A}{1+A}>0
\]
and is a decreasing sequence. Therefore it has a limit
$d_{\infty}>0$. Geometrically, if $\mathcal{C}$ is the boundary of
the pseudo-hyperbolic disk $\Delta=\Delta (i,d_{\infty})$, then
$\mathcal{C}$ intersects the horizontal line $\{\ima z=A \}$ in
one or two points, $q^{+}$ and $q^{-}$, which are the only points
where the sequence $q_{n}$ can accumulate. If $q^{+}$ and $q^{-}$
happen to coincide then that is the limit of  $q_{n}$. Moreover,
if $q^{+}$ and $q^{-}$ are distinct, then let
\[
B=\max_{\zeta \in\mathcal{C}}\ima \zeta >A.
\]
So one can choose $n_{0}$ so that $\ima q_{n}<B$
for $n\geq n_{0}$. Hence the tail $\{q_{n} \}_{n\geq n_{0}}$ cannot
jump from $q^{+}$ to $q^{-}$ because the whole sequence $q_{n}$ stays
in the complement of $\De$. Therefore, we have shown that $q_{n}$
always has a finite limit, which we call $q_{\infty}$. For future use
we note here  that
$q_{\infty}=b_{\infty }+iA$ where
\begin{equation}\label{eq:b}
b_{\infty}=\lim_{n\rightarrow \infty}\frac{x_{n+1}-x_{n}}{y_{n}}.
\end{equation}

Now, since $\tau_{n}$ is a translation followed by a dilation the slope of
the straight segment $[z_{n},z_{n+1}]$ is the same as the slope of
$[i,q_{n}]$, hence we get
\[
\Arg (z_{n+1}-z_{n})\longrightarrow \Arg (q_{\infty}-i).
\]
Fix $\epsilon >0$ and consider the angular sector
\[
S_{\epsilon}=\{z\in\Hh :|\Arg z -\Arg (q_{\infty}-i)|<\epsilon \}.
\]
Then, there exists $n_{0}=n_{0} (\epsilon)$ such that for $n\geq n_{0}$,
$z_{n}$ belongs to the shifted sector
$z_{n_{0}}+S_{\epsilon}$. Letting $n$ tend to infinity we get
\[
 \Arg (q_{\infty}-i)-\epsilon \leq \liminf_{n\rightarrow \infty}\Arg
z_{n}\leq \limsup_{n\rightarrow \infty}\Arg z_{n}\leq \Arg
(q_{\infty}-i)+\epsilon.
\]
This is geometrically clear but can also be seen from the formula
\[
\Arg (z)=\Arg (z-z_{n_{0}})+\Arg \left(1+\frac{z_{n_{0}}}{z-z_{n_{0}}}
\right)=\Arg (z-z_{n_{0}})+o (1)
\]
as $z$ tends to infinity.
Finally, since $\epsilon$ was arbitrary we obtain
\[
\lim_{n\rightarrow \infty}\Arg z_{n}=\Arg (q_{\infty}-i).
\]

\subsection{
%\noindent
{\bf Property (2) (b):} {\em For every angle $\theta \in (0,\pi)$
one can find an orbit $z_{n}$ of $\Phi$ such that $\theta
(z_{0})=\lim_{n\rightarrow \infty}\Arg z_{n}$ is equal to
$\theta$.} }\label{ssec:semiconf}

This property is best established by constructing a conjugation
(change of variables) in the spirit of K\oe nigs' Theorem in the
elliptic case, see (\ref{eq:koenigs}). The existence of such a
conjugation is by itself very interesting, and, after the original
work of Valiron, many others authors deal with such a problem.
Valiron finds a map $\sigma $ such that
\begin{equation}\label{eq:valironsi}
\sigma \circ \Phi = A\sigma
\end{equation}
by showing that the normalized sequence of iterates $\Phi_{n}
(z)/|\Phi_{n} (z_{0})|$ converges uniformly to it. We will use a
slightly different normalization suggested by Pommerenke in
\cite{pom} which has been found useful in other situations, namely
for the parabolic case \cite{pom}, and for backward iterates
\cite{finn} and \cite{iber}. We also recall the work by
Cowen~\cite{Co}, where a different approach, based on the 
uniformization theorem, is used.

The strategy is to renormalize the iterates of $\Phi$ using the
automorphisms $\tau_{n}$ introduced in (\ref{eq:taun}), i.e., choose
an orbit $z_{n}=x_{n}+iy_{n}$ and then study the convergence of the
sequence $\sigma_{n}=\tau_{n}\circ\Phi_{n}$. Observe that $\sigma_{n}
(z_{0})=i$ for all $n=1,2,3,\dotsc$, and since $\sigma_{n}
(z_{1})=q_{n}$ as in (\ref{eq:qn}) we also have
\[
\lim_{n\rightarrow \infty}\sigma_{n} (z_{1})=b_{\infty}+iA.
\]
In particular, every normal sublimit of $\sigma_{n}$ is a non-constant
analytic function.

We first claim that $d (\sigma_{n},\sigma_{n+1})$ tends to $0$ as $n$
tends to infinity. By Schwarz's Lemma,
\[
d(\sigma_n(z), i)=d (\sigma_{n} (z),\sigma_{n} (z_{0}))\leq d
(z,i).
\]
So $\si_{n} (z)$ stays in a compact subset of $\Hh$ and since
\[
\sigma_{n+1} (z)= (\tau_{n+1}\circ \Phi \circ \tau_{n}^{-1})
(\sigma_{n} (z)),
\]
it will be enough to show that the sequence $\psi_{n}=\tau_{n+1}\circ \Phi
\circ \tau_{n}^{-1}$ converges uniformly on compact subsets of $\Hh $
to the identity. Write
\[
\psi_{n}(z)=\frac{\Phi (x_{n}+zy_{n})-x_{n+1}}{y_{n+1}}=
z\frac{y_{n}}{y_{n+1}}\frac{\Phi
(x_{n}+zy_{n})}{x_{n}+zy_{n}}+
\frac{x_{n}}{y_{n}}\frac{y_{n}}{y_{n+1}}\frac{\Phi
(x_{n}+zy_{n})}{x_{n}+zy_{n}}-\frac{x_{n+1}}{y_{n+1}}.
\]
For fixed $z$ the sequence $x_{n}+zy_{n}$ tends to infinity
non-tangentially, so we can apply Julia-Carath\'eodory's Theorem
\ref{thm:jc}, and using
the fact that $x_{n}/y_{n}=\cot\Arg z_{n}$ has a limt, we
obtain that $\psi_{n} (z)$ tends to $z$.

This implies that if $\sigma_{N}$ is a subsequence converging to a
normal sublimit $\sigma$, then $\sigma_{N+1}$ will tend to $\si$ as well.
Therefore, since
\[
\sigma_{n}\circ \Phi = (\tau_{n}\circ \tau_{n+1}^{-1})\circ \sigma_{n+1}
\]
for all $n$  and since
\[
\tau_{n}\circ \tau_{n+1}^{-1}
(z)=\frac{x_{n+1}-x_{n}}{y_{n}}+z\frac{y_{n+1}}{y_{n}}\rightarrow
b_{\infty }+Az
\]
by (\ref{eq:b}) and (\ref{eq:yn}), we obtain that every sublimit
$\sigma$ must satisfy the functional equation
\begin{equation}\label{eq:functeq}
\sigma \circ \Phi =A\sigma +b_{\infty}.
\end{equation}

Finally, since
\[
d (i,\sigma_{n} (z))=d (\tau_{n}\circ \Phi_{n} (z_{0}),\tau_{n}\circ
\Phi_{n} (z))= d (\Phi_{n} (z_{0}),\Phi_{n} (z) )
\]
is a decreasing sequence, it must converge to $d (i,\sigma
(z))$. Hence any other normal sublimit $\tilde{\sigma}$ must satisfy $d
(i,\tilde{\sigma}(z))=d (i,\sigma(z))$ for all $z\in \Hh $,
i.e., $\tilde{\sigma}$ can
only differ from $\sigma$ by an automorphism of $\Hh$ which fixes $i$.
However, writing $T (z)=Az+b_{\infty}$, equation (\ref{eq:functeq})
can be iterated to $\sigma \circ \Phi_{n} =T_{n}\circ \sigma$, hence
we get that $\sigma (z_{n})=T_{n} (i)$ which is a sequence tending to infinity.
In particular, $\tilde{\si}$ can only differ from $\sigma$ by  an
automorphism of $\Hh$ which fixes $i$ and infinity, but this can only
be the identity.

In conclusion, we have shown that given an orbit $z_{n}$ of $\Phi $
one can renormalize the iterates of $\Phi$
with some automorphisms $\tau_{n}$ of $\Hh$ built from $z_{n}$ so that
$\tau_{n}\circ \Phi_{n}$ converges uniformly on compact subsets of
$\Hh$ to a function $\sigma$ which satisfies the functional equation
\begin{equation}\label{eq:pommersi}
\si \circ \Phi = T\circ \sigma = A \si +b_{\infty}
\end{equation}
where $T (z)=Az+b_{\infty}$ and
$b_{\infty}$ is a real number depending continuously on $z_{0}$.
In fact, if $\theta (z_{0})=\lim_{n\rightarrow \infty}\Arg \Phi_{n}
(z_{0})$, see Property 2 (a) above, then
\begin{equation}\label{eq:cot}
b_{\infty}= (A-1)\cot
(\theta (z_{0})).
\end{equation}

Writing $\hat{\sigma}=\sigma +b_{\infty}/ (A-1)$, one sees that
$\hat{\sigma}$ satisfies (\ref{eq:valironsi}), and a computation
using (\ref{eq:cot}) shows that actually $\hat{\si}$ and Valiron's
conjugation are the same function.  Yet, one may ask: how many
solutions do (\ref{eq:valironsi}) and (\ref{eq:pommersi}) have?
Also, is it possible to choose $z_{0}$ so that in
(\ref{eq:pommersi}) the coefficient $b_{\infty}$ becomes $0$?
Namely, we still haven't established Property 2~(b).

\noindent {\bf Semi-conformality of $\sigma$}. All the previous
questions can be answered if we can show that the conjugating map
$\sigma$ that we have found in (\ref{eq:pommersi}) has the
property of being {\em semi-conformal} at infinity. Without loss
of generality we can work with $\hat{\sigma}$ instead of $\sigma$.
Thus we want to show that $\ntlim_{z\rightarrow
\infty}\hat{\sigma} (z)=\infty$ and that
\begin{equation}\label{eq:semiconf}
\ntlim_{z\rightarrow \infty}\Arg \frac{\hat{\sigma} (z)}{z}=0.
\end{equation}
To this end we introduce the functions
\[
g_{n}=A^{-n} \hat{\si}\circ \tau_{n}^{-1}-\frac{b_{\infty}}{A-1}
\]
which are self-maps of $\Hh$ (we follow the same argument as in
Section 2 of \cite{finn}). Notice that, since
$\hat{\sigma}(z_0)=\sigma(z_0)+b_\infty/(A-1)=i+b_\infty/(A-1)$,
then
\[
g_{n} (i)=A^{-n}\hat{\si}(z_{n})-\frac{b_{\infty}}{A-1}=i.
\]
Also if $q_{n}$ is defined as in (\ref{eq:qn}), then
\[
g_{n} (q_{n})=A\hat{\si} (z_{0})=Ai+b_{\infty}=q_{\infty}.
\]
So any normal sublimit $g$ of the sequence $g_{n}$ must fix $i$
and also $q_{\infty}$ (since $q_{n}\rightarrow q_{\infty}$). Thus
by Schwarz's lemma, $g$ is the identity, i.e., $g_{n}
(z)\rightarrow z$ uniformly on compact subsets of  $\Hh$. Now let
$K$ be a compact subset of $\Hh$. As $n\rightarrow \infty $,
\[
d \left(A^{n} \left(z+\frac{b_{\infty}}{A-1} \right), \hat{\si}
(x_{n}+zy_{n}) \right) =d (z,g_{n} (z))\rightarrow 0
\]
uniformly for $z\in K$, hence
\[
\Arg \hat{\si} (x_{n}+zy_{n})-\Arg\left(z+\frac{b_{\infty}}{A-1}
\right) \rightarrow 0.
\]
But by (\ref{eq:cot}) we also have
\[
\Arg (x_{n}+y_{n}z)=\Arg (z+\frac{x_{n}}{y_{n}})\rightarrow \Arg
(z+\cot\theta (z_{0}))=\Arg \left(z+\frac{b_{\infty}}{A-1}\right).
\]
Hence,
\[
\Arg \hat{\si} (x_{n}+zy_{n})-\Arg (x_{n}+y_{n}z)\rightarrow 0.
\]
By choosing $K$ to be a hyperbolic disk of larger and larger
radius we see that the union of the sets $x_{n}+y_{n}K$ eventually
covers sectors of larger and larger opening. We have proved the
semi-conformality of $\hat{\si}$, and thus of $\sigma$.

\smallskip

Now that we know that $\hat{\sigma}$ is semi-conformal, iterating
(\ref{eq:valironsi}), which is satisfied by $\hat{\sigma}$, we get
$\hat{\sigma } \circ \Phi_{n}=A^{n}\hat{\sigma}$ and evaluating at
$z_{0}$ we obtain $\hat{\sigma} (z_{n})=A^{n}\hat{\si}(z_{0})$.
Applying (\ref{eq:semiconf}) to $z_{n}$ we see that $\Arg
A^{n}\hat{\sigma}(z_{0})-\Arg z_{n}$ tends to zero. In other
words, $\theta (z_{0})=\Arg\hat{\sigma}(z_{0})$. It then remains
to show that by varying $z_{0}$, $\Arg\hat{\sigma}(z_{0})$ takes
on every value in $(0,\pi)$. This follows at once from the
semiconformality as well and it is explained in Lemma
\ref{lem:semiconf} below. For a proof of this lemma see Section 5
of \cite{finn}, and also \cite{Co}.
\begin{lemma}\label{lem:semiconf}
Suppose $\hat{\si}$ is an analytic self-map of $\Hh$ which has
non-tangential limit infinity at infinity and is semi-conformal, i.e.,
(\ref{eq:semiconf}) holds for $\hat{\si}$. Then there is a simply-connected
region $\Omega$ in $\Hh$ with an {\em inner-tangent at infinity}, i.e. for
every $\alpha \in(0,\pi/2)$ there is $R>0$ so that
\[
\{|\Arg z-\pi /2|<\alpha;\ |z|>R \}\subset\Omega,
\]
with the property that $\hat{\si}$ restricted to $\Om$ is one-to-one and
$\hat{\si} (\Om)$ also has an inner tangent at infinity.
\end{lemma}

The previous lemma in particular allows to select a
simply-connected region $\Omega\subset \Hh$, called a {\em
fundamental set} for $\phi$, so that
\begin{enumerate}
\item The map $\phi$ is one-to-one on $\Omega$.
\item The set $\Omega$ is {\em fundamental} for $\phi$, in the
sense that $\phi(\Omega)\subseteq \Omega$ and for any compact
subset $K\subset\subset \Hh$ there exists $N=N(K)$ so that
$\phi_n(K)\subset \Omega$ for all $n>N$.
\item The map $\hat{\sigma}$ is one-to-one on $\Omega$.
\item The set $\hat{\sigma}(\Omega)$ is fundamental for the
hyperbolic automorphism $\zeta \mapsto A \zeta$.
\end{enumerate}

The explicit knowledge of the set $\Omega$ (and the intertwining
map $\sigma$) coincides with the knowledge of the analytic and
dynamical properties of $\phi$. One could say that the dynamical
properties of $\phi$ are read by means of the geometrical
properties of the couple $(\sigma, \Omega)$. For instance, $\phi$
is one-to-one on $\Hh$ if and only if $\Omega=\Hh$ if and only if
$\sigma$ is one-to-one on $\Hh$.

To go back to our questions, we are left to deal with the
uniqueness of the map $\sigma$. We have

\begin{proposition}[Uniqueness of conjugation]\label{prop:unique}
Suppose $\sigma $ is an analytic self-map of $\Hh$ which satisfies
the functional equation (\ref{eq:valironsi}). Then $\sigma$ has
non-tangential limit $\infty$ at $\infty$, it is semi-conformal at
$\infty$ ({\em i.e.} \eqref{eq:semiconf} holds for $\sigma$).
Moreover, every other self-map of $\Hh$ satisfying
(\ref{eq:valironsi})  is a positive constant multiple of $\sigma$.
\end{proposition}

Oddly enough, even if everyone would swear that all the solution
built by Valiron~\cite{va1}, Pommerenke~\cite{pom},
Cowen~\cite{Co} and Bourdon-Shapiro~\cite{bs} coincide, it seems
that no one proved  this explicitly. In case the map $\sigma$ is
known to fix $\infty$ as non-tangential limit and to be
semi-conformal at $\infty$, the proof can be done directly (see
the proof of Theorem~1.2 of \cite{finn}). Here we present a
different proof which is based on the existence of an intertwining
map semi-conformal at $\infty$ and a theorem on the commutator of
hyperbolic automorphisms due to Heins~\cite{He}.

\noindent{\bf The Proof of Proposition~\ref{prop:unique}}. Let us
first define the following sets of holomorphic mappings:
\[
{\mathcal C}_A:=\{F : \Hh \to \C \ \hbox{holomorphic}\ | F (Az)=A
F(z)\ \  \forall z\in \Hh\},
\]
\[
\mathcal S:=\{ \sigma : \Hh \to \C \ \hbox{holomorphic}\ | \sigma
\circ \phi = A \sigma \}.
\]
The set $\mathcal C_A$ is thus formed by holomorphic maps which
commute with the linear fractional map (hyperbolic automorphism of
$\Hh$) $\zeta \mapsto A \zeta$; while the set $\mathcal S$ is made
of all solutions of the functional equation~\eqref{eq:valironsi}.
Notice that for the moment we are not restricting ourselves to
{\sl self-maps} of $\Hh$. The two sets are essentially the same as
the following lemma shows (see also Lemma~4 in~\cite{Co}). As a
matter of notation, we let $\sigma_V$ be the Valiron intertwining
mapping constructed before.

\begin{lemma}\label{corr}
There is a one-to-one correspondence between $\mathcal S$ and
$\mathcal C_A$ given by:
\[
\mathcal C_A \ni F \mapsto F \circ \sigma_V \in \mathcal S.
\]
\end{lemma}

\begin{proof}
Let $F \in \mathcal C_A$. Let us denote by $\Phi(z)=Az$. Then
\[
(F \circ \sigma_V) \circ \phi= F \circ (\sigma_V \circ \phi)= F
\circ \Phi \circ \sigma_V = \Phi \circ (F \circ \sigma_V).
\]
On the other hand if $\sigma \in \mathcal S$, since $\sigma_V$ is
univalent on $\Omega$ (the fundamental set constructed before),
one can define a holomorphic map $\tilde{F}$ on $\sigma_V(\Omega)$
by
\[
\tilde{F}(\sigma_V(x)):=\sigma \circ \sigma_V^{-1}(x).
\]
Since $A(\sigma_V(\Omega))\subseteq \sigma_V(\Omega)$, on
$\sigma_V(\Omega)$ we have
\begin{equation}\label{primo}
\tilde{F} \circ \Phi=\sigma\circ\sigma_V^{-1}\circ \Phi = \sigma
\circ \phi \circ \sigma_V^{-1}=\Phi \circ \sigma \circ
\sigma_V^{-1} = \Phi \circ \tilde{F}.
\end{equation}
Then one can extend $\tilde{F}$ to all of $\Hh$ as follows:
\[
F(z)=A^{-n}\tilde{F}(A^n z) \quad \hbox{for $z \in \Hh$ and $n\in
\N$ such that $A^n z \in \sigma_C(V)$}.
\]
The map $F$ is well defined, {\sl i.e.}, it is independent of
$n\in \N$ by~\eqref{primo}. Moreover $F \in \mathcal C_A$ and
$\sigma=F\circ \sigma_V$.
\end{proof}

Now we can complete the proof of Proposition~\ref{prop:unique} as
follows. Let $\sigma\in\mathcal S$ be such that
$\sigma(\Hh)\subseteq \Hh$. From Lemma~\ref{corr} it follows that
$\sigma=F \circ \sigma_C$ for some $F:\Hh\to\C$ such that $F(A
w)=A F(w)$.  If $F(\Hh)\subseteq \Hh$, by a theorem of
Heins~\cite{He}, we must have $F(w)=\mu w$ for some $\mu \in \R^+$ and
therefore $\sigma=\mu \sigma_V$, which in particular proves that
$\sigma$ has fixed point $\infty$ and it is semi-conformal at
$\infty$. We are thus left to prove that if
$F(\sigma_V(\Hh))\subseteq \Hh$ then actually $F(\Hh)\subseteq
\Hh$. Assume this is not the case. Then there exists $w_0\in \Hh$
such that $\ima F(w_0)\leq 0$. Since $\sigma_V(\Hh)$ is
fundamental for $w\mapsto A w$, it follows that there exists $n\in
\N$ such that $A^n w_0\in \sigma_V(\Hh)$. But then
\[
\ima  F(A^n w_0)=\ima A^n F(w_0) \leq 0,
\]
meaning that $F \circ \sigma_V(\Hh) \not\subset \Hh$ against our
hypothesis.

\begin{remark}
More generally, arguing as in Proposition~4 of~\cite{Co} one can
prove that for any $\sigma \in \mathcal S$ (no restriction on the
image $\sigma(\Hh)$) there exists a holomorphic map $g:\{\zeta \in
\C : |\log|\zeta||<\pi^2/\log A\}\to \C$ such that $\sigma$ is
given by $w\mapsto \sigma_V(w)\cdot g(\exp(2\pi i \log
\sigma_V(w)/\log A))$. Thus Proposition~\ref{prop:unique} says
that if $\sigma(\Hh)\subseteq \Hh$ then $g$ is a real positive
constant.
\end{remark}

\subsection{{\bf Property 3:} {\em The ratios $z_{n}/A^{n}$ do not
always converge.}}\label{ssec:conf} This property is connected to the
conformality at infinity of Valiron's conjugation. In fact,
let $\sigma$ be the limit of
$\tau_{n}\circ \Phi_{n}$, and without loss of generality assume
that $b_{\infty}=0$ so that $\sigma$ satisfies (\ref{eq:valironsi}).
Let $\alpha =\inf_{z\in \Hh}\ima \sigma
(z)/\ima z$. There are two possibilities: either $\alpha $ is $0$ or
it is positive. In either case, since $z_{n}$ approaches infinity
non-tangentially,  Julia-Carath\'eodory's Theorem \ref{thm:jc}
applies, so that
\[
\frac{\sigma (z_{n})}{z_{n}}=\frac{A^{n}}{z_{n}}\si (z_{0})\rightarrow \alpha
\]
When $\alpha >0$ it is costumary to say that $\sigma$ has a finite
angular derivative at infinity. Valiron gives a couple of
necessary and sufficient conditions for this to happen, which are
quite tautological. Bourdon and Shapiro \cite{bs} show that if
$\Phi$ extends analytically near infinity then $\alpha >0$.
Arguing as in~\cite{bg} one can state the Bourdon-Shapiro theorem
as follows:
\begin{theorem}\label{thm:bs}
Suppose $\Phi$ is an analytic self-map of $\Hh$ such that $\Phi
(z)=Az+\Ga (z)$, with $A>1$, and there exist $M,\epsilon>0$ such
that  $|\Ga (z)|\leq M|z|^{1-\epsilon}$ for all $z\in \Hh$. Then
the conjugating map $\sigma$ has a finite angular derivative at
infinity.
\end{theorem}

The question of the convergence of the ratio $z_{n}/A^{n}$ is
strictly related to that of the existence of fixed points for
intertwining mappings $\sigma:\Hh\to \Hh$. Indeed, assume that
$\sigma_V$ has finite angular derivative at infinity, say
$\alpha>0$. Then $\sigma_\lambda:=\lambda\sigma_V$ for
$\lambda>1/\alpha$ is a holomorphic self-map of $\Hh$ such that it
has non-tangential limit $\infty$ at $\infty$, and
$\sigma_\lambda'(\infty)=\lambda \alpha >1$. Therefore $\infty$ is
the Denjoy-Wolff point of $\sigma_\lambda$ for all
$\lambda>1/\alpha$. In particular $\sigma_\lambda$ has no fixed
points in $\Hh$. Therefore, {\sl the ratio $z_{n}/A^{n}$ is
convergent if and only if there exists one---and hence infinitely
many---intertwining maps $\sigma$ with Denjoy-Wolff point at
$\infty$}. One is thus forced to study the following curve $T$:
\[
T:\R^+\ni t \mapsto H(t \sigma_V)\in \overline{\Hh}\cup\{\infty\},
\]
where for a holomorphic self-map $f\neq Id$ of $\Hh$, $H(f)$ is
the so-called {\sl Heins map}, defined to be the (unique) fixed
point of $f$ in $\Hh$ if $f$ has fixed points, or the Denjoy-Wolff
point of $f$ in case $f$ has no fixed points in $\Hh$. The map $H$ is
easily seen to be continuous on the subset of the complex Banach
space $H^\infty(\Hh)$ given by functions with range in
$\overline{\Hh}$, and it can be shown that it is holomorphic on
the open set given by functions whose image is relatively compact
in $\Hh$ (see~\cite{br}). Therefore the curve $t\to T(t)$ is a
continuous curve in $\overline{\Hh}$ that can be continuously
extended to $[0,\infty)$ as $T(0)=0$ (the geometric meaning is
that the constant function $z\mapsto 0$ is a solution
of~\eqref{eq:valironsi}). Moreover it is analytic at a point $t_0$
whenever $T(t_0)\in \Hh$. The question on the ratio $z_{n}/A^{n}$
can be stated in terms of $T$ as follows: {\sl the ratio
$z_{n}/A^{n}$ is convergent if and only if the curve $T$ reaches
infinity in a finite time, namely if and only if  there exists
$t_0 \in (0,+\infty)$ such that $T(t_0)=\infty$} (and then
$T(t)=\infty$ for $t>t_0$). The curve $T$ reads the geometrical
properties of $\phi$. For instance it is easy to see that if
$\phi$ is such that $\lim_{z\to p} |\phi(z)|<1$ for all $p \in
\partial \Hh\setminus \{\infty\}$, then $T(t)\in \Hh \cup
\{\infty\}$ for all $t\in (0,+\infty)$, and in particular $T$ is
analytic in its interior. With a slightly more subtle argument on
commuting mappings (using Behan's lemma, see, {\sl e.g.},
\cite{ab}) one can show that if $T(t)=0$ for some $t\in
(0,+\infty)$ then $f$ cannot fix $0$ in the sense of
non-tangential limits. It would be interesting to pursue a
systematic study about the relations between properties of $\phi$
and properties of $T$.

\section{Several complex variables}\label{sec:svc}

We fix $N=2,3,4,\dotsc$ and $\B=\B^{N}=\{z\in \C^N: \|z\|<1 \}$, where
\[
\|z\|^{2}= (z,z) \qquad\mbox{ and }\qquad
(z,w)=\sum_{j=1}^{N}z^{j}\overline{w^{j}}.
\]
Let $\phi$ be a self-map of $\B$. As in the disk case we can say that
$\phi$ is of {\em elliptic type} if it fixes at least one point in
$\B$ (however, now, $\phi$ could fix more than just one-point and not
be the identity). We are interested in the case when $\phi$ has no fixed
points in $\B$. The Denjoy-Wolff Theorem still hold (see \cite{ab}
Theorem 2.2.31), namely, the iterates of $\phi$ converge to one point
on $\bd{\B}$. By conjugating with a unitary map we can assume without
loss of generality that this special point is $e_{1}=(1,0,\dots,0)$.
Once again maps with no fixed points in $\B$ will be divided into
hyperbolic and parabolic type, but before we
can do this we need to introduce a few tools.

\subsection{A special automorphism}\label{ssec:spaut}
For $a\in \B$, we define the projections
\[
P_{a} (z)=\frac{(z,a)}{(a,a)}a\qquad \mbox{ and }\qquad Q_{a}
(z)=z-P_{a} (z).
\]
Then we let
\begin{equation}\label{eq:ga}
\ga _{a} (z)=\frac{P_{a} (z)+s_{a}Q_{a} (z)-a}{1- (z,a)}
\end{equation}
where $s_{a}=\sqrt{1-\|a\|^{2}}$, and so that $\ga_{a} (a)=0$. It is
well-known that $\ga_{a}$ is an automorphism of $\B$.

We define the pseudo-hyperbolic distance between two points $a,b \in
\B$ as
\[
d (a,b)=\|\ga_{a} (b)\|<1.
\]
Schwarz's Lemma (\cite{ab} Thm. 2.2.12) and \cite{ab} Cor. 2.2.2,
imply as in the disk that $d (\phi (a),\phi (b))\leq d (a,b)$. Another
quantity which is decreased by self-maps of the ball is
\begin{equation}\label{eq:qab}
Q (a,b)=\frac{|1- (a,b)|^{2}}{(1-\|a\|^{2}) (1-\|b\|^{2})},
\end{equation}
i.e., $Q (\phi (a),\phi (b))\leq Q (a,b)$ (\cite{ab} Prop 2.2.17)

\subsection{Hyperbolic versus parabolic}\label{ssec:hypvpar}
We will again consider the orbit of the origin
$z_{n}=\phi_{n} (0)$, thus $z_{n}\rightarrow e_{1}$. It follows that
one can extract a subsequence $z_{N}$ with the property that
$\|z_{N+1}\|\geq \|z_{N}\|$. Hence
\[
c =\liminf_{z\rightarrow e_{1}}\frac{1-\|\phi (z)\|}{1-\|z\|}\leq 1
\]
so by Julia's Lemma (\cite{ab} Thm. 2.2.21)
\[
\frac{|1- (\phi (z),e_{1})|^{2}}{1-\|\phi  (z)\|^{2}}
\leq c \frac{|1- (z,e_{1})|^{2}}{1-\|z\|^{2}},
\]
and in particular,
\begin{equation}\label{eq:jlball}
\frac{|1- z_{n+1}^1|^{2}}{1-\|z_{n+1}\|^{2}}
\leq c \frac{|1- z_n^1 |^{2}}{1-\|z_{n}\|^{2}}.
\end{equation}
We say $\phi$ is of {\em hyperbolic type} if
$c <1$, and of {\em parabolic type} if $c=1$. The quantity $c=c
(\phi)$ is called the coefficient of dilatation of $\phi$.
It is clear that
\begin{equation}\label{eq:dilcoeff}
c (\phi_n)=[c (\phi)]^{n}.
\end{equation}

In the sequel we will assume that $\phi$ is a self-map of the ball of
hyperbolic type. First we describe the automorphisms of hyperbolic type.

\subsection{Automorphisms of the ball of hyperbolic type}\label{ssec:authyp}
As in the one-dimensional case it is best to move to an ``upper
half-plane'' model. It turns out that $\B$ is biholomorphic to the
domain
\[
\Hh^{N}=\{w= (w^{1},w^{\prime})\in \C^N:\ \ima w^{1}>\|w^{\prime}\|^{2} \}
\]
via a map very similar to the classical Caley transform.  Given $a\in
\Hh^{N}$ with $\ima a^{1}-\|a^{\prime}\|^{2}>1$ there is an
automorphism of hyperbolic type $\Psi_a$ which sends the point $\iota
=(i,0^{\prime })$ to $a$.  We first build the inverse of such mapping.
Consider the translation (we refer to \cite{ab} p.~155 for these
automorphisms of $\Hh^{N}$.)
\begin{equation}\label{eq:transl}
h_{b} (w)= (w^{1}+b^{1}+2i\langle w^{\prime },b^{\prime
}\rangle,w^{\prime }+b^{\prime })
\end{equation}
where $b= (-\rea a^{1}+i\|a^{\prime } \|^{2},-a^{\prime })\in
\bd\Hh^{N}$. Then
\[
h_{b} (a)= (i (\ima a^{1} -\|a^{\prime}\|^{2}),0^{\prime }).
\]
Now consider the non-isotropic dilation
\begin{equation}\label{eq:dilation}
\delta _{A} (w)=\left(Aw^{1},\sqrt{A}w^{\prime } \right)
\end{equation}
where $A=\ima a^{1} -\|a^{\prime }\|^{2}$.
The automorphism $\Phi_{a}=\delta _{1/A}\circ h_{b}$ sends $a$
to $\iota$,
\[
\Phi_{a} (w)= \left(\frac{w^{1}-\rea a^{1}+i\|a^{\prime }
\|^{2}-2i\langle w^{\prime },a^{\prime }\rangle}{A},
\frac{w^{\prime }-a^{\prime }}{\sqrt{A}} \right).
\]
The inverse is
\begin{equation}\label{eq:psia}
\Psi_a (z)=\left(Az^{1}+\rea a^{1}+i\|a^{\prime }\|^{2}+2i\sqrt{A}\langle
z^{\prime },a^{\prime }\rangle, \sqrt{A}z^{\prime }+a^{\prime }\right).
\end{equation}
More generally, given a unitary transformation $U$ of $\C^{N-1}$ one
can consider the automorphism
\begin{equation}\label{eq:hypaut}
\Psi(z)=\left(Az^{1}+\rea a^{1}+i\|a^{\prime }\|^{2}+2i\sqrt{A}\langle
U (z^{\prime }),a^{\prime }\rangle, \sqrt{A}U (z^{\prime })+a^{\prime }\right).
\end{equation}
Varying $a$ with $A=\ima a^{1}-\|a^{\prime}\|^{2}>1$ and $U$ as above,
the automorphisms $\Psi$ describe all possible hyperbolic
automorphisms of $\Hh^{N}$ with infinity as attracting fixed point.

\subsection{Self-maps of the ball of hyperbolic
type}\label{ssec:smball} Let $\Phi $ be a holomorphic self-map of
$\Hh^{N} $ without fixed points in $\Hh^{N} $, such that its
Denjoy-Wolff point is $\infty$, and of hyperbolic type.

Following the
lead of the one-dimensional case the following open problem arises:
\begin{openpb}\label{oppb:valscv}
How closely are the orbits of $\Phi$ trying to imitate the behavior of
the orbits of a corresponding hyperbolic automorphism $\Psi$ as in
(\ref{eq:hypaut})?
\end{openpb}

The automorphism $\Psi$ in (\ref{eq:hypaut}) fixes exactly two points:
infinity and the point $c\in \bd{\Hh^{N}}$. To see this first solve
\[
\sqrt{A}U (c^{\prime})+a^{\prime}=c^{\prime}.
\]
Taking $U^{-1}$ and dividing by $\sqrt{A}$, one gets
\[
(I-\frac{1}{\sqrt{A}}U^{-1}) (c^{\prime})=-\frac{1}{\sqrt{A}}U^{-1}
(a^{\prime})
\]
which is invertible. Now solve
\[
Ac^{1}+\rea a^{1}+i\|a^{\prime }\|^{2}+2i\sqrt{A}\langle
U (c^{\prime }),a^{\prime }\rangle=c^{1}
\]
using the fact that $\sqrt{A}U (c^{\prime})= (I-
(\sqrt{A}U)^{-1})^{-1} (a^{\prime})$.

Therefore using an appropriate translation as in (\ref{eq:transl}),
the map $\Psi$ can be conjugated to an automorphism whose fixed points
are $0$ and $\infty$, i.e. to
\begin{equation}\label{eq:psio}
\Psi_0 (z)= (Az^{1},\sqrt{A}U(z^{\prime})).
\end{equation}
Moreover, by linear algebra $\Psi_0$ can be further conjugated via a
unitary matrix so that $U$ becomes diagonal.

Our open problem can be rephrased as
\begin{openpb}\label{oppb:valscv2}
Given a holomorphic self-map $\Phi$ of
$\Hh^{N}$ without fixed points in $\Hh^{N}$, such that its
Denjoy-Wolff point is $\infty$, and which is of hyperbolic type with
dilation coefficient $A>1$, does there exist a unitary trasformation
$U$ of $\C^{N-1}$ and a conjugation $\sigma$ (also a self-map of
$\Hh^{N}$) such that
\[
\sigma \circ \Phi = \Psi_0\circ \sigma
\]
where $\Psi_0$ is as in (\ref{eq:psio}), and so that $\sigma$ has some
degree of regularity at infinity to be determined (something along the
lines of semi-conformality)?
\end{openpb}

In a recent preprint \cite{bg} the result  of Bourdon and Shapiro
has been generalized to several complex variable, {\sl i.e.},
Valiron's conjugation is established under some smoothness
assumptions at infinity for $\Phi$.

What we can show is a partial answer to Open Problem \ref{oppb:valscv}
which resembles
Lemma \ref{lem:ntapproach} in the one-dimensional case. Namely, we can
show that the orbits of $\Phi$ always remain in a  Kor\'anyi
approach region at infinity, see definition below
(the fact that the orbits of $\Psi_{0}$ remain in a Kor\'anyi
approach region at infinity can easily be verified).

\subsection{Kor\'anyi approach of the orbits}\label{ssec:kapproach}
Back in the ball setting, let $\phi $ be a holomorphic self-map of
$\B$ without fixed points in $\B$, such that its
Denjoy-Wolff point is $e_{1}$, and which is of hyperbolic type with
dilation coefficient $c<1$.

Given a parameter $M>0$ the Kor\'anyi regions at $e_{1}$ of amplitude
$M$ are the sets
\[
K (R)=\left\{z\in \B : \frac{|1-z^{1}|}{1-\|z\|^{2}}<R \right\}.
\]
We need a preliminary result.
\begin{claim}\label{cl:Kapproach}
If $c<3-\sqrt{8}$,
then the orbit $z_{n}=\phi_n (0)$ tends to $e_{1}$ while staying in a
Kor\'anyi approach region, i.e.,
\[
L_{n}=\frac{|1- (z_{n},e_{1})|}{1-\|z_{n}\|^{2}}
=\frac{|1- z_n^1 |}{1-\|z_{n}\|^{2}}\leq M<\infty
\]
for some constant $M<\infty$.
\end{claim}
Assuming Claim \ref{cl:Kapproach} for the moment, we show the
Kor\'anyi approach of the orbit $z_{n}=\phi_n(0)$.
Using (\ref{eq:dilcoeff}), we can
find an integer $N$ large enough so that $c^{N}<3-\sqrt{8}$, and Claim
\ref{cl:Kapproach}  implies that $z_{kN}$, $k=1,2,3\dotsb$, stays in a
Kor\'anyi region. However, for $j=1,\dotsc ,N-1$, $d
(z_{kN+j},z_{kN})\leq d (0,z_{j})$, by Schwarz's Lemma. Hence, since
the hyperbolic neighborhood of a Kor\'anyi region is still a Kor\'anyi
region, we find that the whole orbit $z_{n}$ remains in a Kor\'anyi region.
Moreover, by the same argument, any orbit $\phi_n (z_{0})$ has the
same property.

\begin{proof}[Proof Claim \ref{cl:Kapproach}:]
We rewrite (\ref{eq:jlball}) as
\begin{equation}\label{eq:sn}
S_{n}=\left|\frac{1- z_{n+1}^1}{1- z_n^1 }\right|\leq
c\frac{L_{n}}{L_{n+1}}.
\end{equation}
Recalling the definition and monotonicity property of $Q (a,b)$ given in
(\ref{eq:qab}), we see that $Q(z_{n},z_{n+1})$ is decreasing and
thus
\begin{equation}\label{eq:wge}
Q(z_{n},z_{n+1})=\frac{|1- (z_{n},z_{n+1})|^{2}}{(1-\|z_{n}\|^{2})
(1-\|z_{n+1}\|^{2})} \leq Q
(z_{0},z_{1})=\frac{1}{1-\|z_{1}\|^{2}}<\infty.
\end{equation}
Notice that
\[
1- (z_{n},z_{n+1})= (e_{1}-z_{n},e_{1})+ (e_{1},e_{1}-z_{n+1})-
(e_{1}-z_{n},e_{1}-z_{n+1}).
\]
Therefore,
\[
|1- (z_{n},z_{n+1})|\geq | (1-z_n^1 )+ (1-
\overline{z_{n+1}^1})|-\|e_{1}-z_{n}\|\|e_{1}-z_{n+1}\|.
\]
Expanding the square, we have
\[
\frac{\|e_{1}-z_{n}\|^{2}}{1-\|z_{n}\|^{2}}=2\frac{1-\rea z_n^1
}{1-\|z_{n}\|^{2}} -1\leq 2L_{n}.
\]
So, after a square root and the triangle inequality, (\ref{eq:wge}) becomes
\begin{equation}\label{eq:final}
\frac{| (1-z_n^1 )+ (1-
\overline{z_{n+1}^1})|}{\sqrt{1-\|z_{n}\|^{2}}
\sqrt{1-\|z_{n+1}\|^{2}}} -2\sqrt{L_{n}L_{n+1}}\leq \sqrt{Q
(z_{0},z_{1})}.
\end{equation}
Now suppose that
\[
\limsup_{n\rightarrow \infty}L_{n}=+\infty.
\]
Then one can find a subsequence $L_{N}$ such that $L_{N}\leq L_{N+1}$
and $L_{N}\rightarrow +\infty$. By (\ref{eq:sn}),
$\limsup_{N\rightarrow \infty} S_{N}\leq  c$. On the other hand,
dividing by $\sqrt{L_{N}L_{N+1}}$ and letting $N$ tend to infinity in
(\ref{eq:final}), we also have
\[
\limsup_{N\rightarrow \infty}\frac{| (1-z_n^1 )+ (1-
\overline{z_{n+1}^1})|}{\sqrt{|1-z_{N}^{1}|}
\sqrt{|1-z_{N+1}^{1}|}}\leq 2.
\]
Squaring both sides and reorganizing
\[
\limsup_{N\rightarrow \infty}|1-S_{N}||1-\frac{1}{S_{N}}|\leq 4.
\]
So if $S$ is a sublimit of $S_{N}$ it must satisfy $S\leq  c$ and
\[
(1-S)^{2}\leq 4S,
\]
{\sl i.e.}, $3-\sqrt{8}\leq S\leq 3+\sqrt{8}$. In particular, if
$c$ happens to be less than $3-\sqrt{8}$, then no sublimit of
$S_{N}$ can exists and therefore $L_{n}$ remains bounded.
\end{proof}

\subsection{Conclusion}
The problem one encounters after this claim is established is that one
would like to use the Julia-Carath\'eodory Theorem for self-maps of
the ball. Such result exists, see \cite{ab} Theorem (2.2.29), however
in order to use it one would need a much more restrictive approach for
the orbit of $\phi$: ``special and restricted''. Of course, even the
orbits of the automorphism $\Psi_0$ do not have this property in
general, but there is always one orbit that does. Our hope is to be
able to produce at least one orbit of $\phi$ that has a special and restricted
approach and then renormalize $\phi$ using this orbit.

There is a different approach which seems to bypass the unitary matrix $U$
of Open Problem \ref{oppb:valscv2}.

\begin{openpb}\label{oppb:valnto1}
Given a holomorphic self-map $\Phi$ of
$\Hh^{N}$ without fixed points in $\Hh^{N}$, such that its
Denjoy-Wolff point is $\infty$, and which is of hyperbolic type with
dilation coefficient $A>1$, does there exist
a conjugation $\eta:\Hh^{N}\longrightarrow \Hh $ such that
\[
\eta \circ \Phi = A\eta,
\]
and so that $\sigma$ has some
degree of regularity at infinity to be determined (something along the
lines of semi-conformality)?
\end{openpb}
Of course if one can find $\sigma$ which solves Open Problem \ref{oppb:valscv2}
then $\eta =\pi^{1}\circ \si$, where $\pi^{1}
(z^{1},z^{\prime})=z^{1}$, will solve Open Problem \ref{oppb:valnto1}.

\end{document}